\documentclass[11pt,fleqn]{article}
\usepackage{stmaryrd}
\usepackage{amssymb}
\usepackage{bbm}
\usepackage{mathrsfs}
\usepackage{amsfonts}
\usepackage{amsmath,epsf,eufrak}
\begin{document}
\renewcommand{\theequation}{\arabic {section}.\arabic {equation}}
\newcommand{\bi}{\begin{equation}}
\newcommand{\ei}{\end{equation}}
\date{}
\title{  Proof of the Sendov conjecture for polynomials of degree nine
\author{Zaizhao Meng}}
\maketitle \baselineskip 24pt
 \begin{center}
 {\bf Abstract}
 In this paper, we prove the Sendov conjecture for polynomials of degree nine.
 We use a new idea to obtain new upper bound for the $\sigma-$sum to zeros of the polynomial.
 \end{center}
\section{Introduction}
\setcounter{equation}{0}

We continue the work of [5].
Let $ \mathcal{P}_{n}$ denote the set of all monic polynomials of degree $n(\geq 2)$ of the form
$$p(z)=\prod\limits_{k=1}^{n}(z-z_{k}),\ |z_{k}|\leq 1(k=1,\cdots,n)$$
with\\
$p^{\prime}(z)=n\prod\limits_{j=1}^{n-1}(z-\zeta_{j}),\ \ |\zeta_{j}|\leq 1(j=1,\cdots,n-1).$\\
Write $I(z_{k})=\min\limits_{1\leq j\leq n-1}|z_{k}-\zeta_{j}|,\ \ I(p)=\max\limits_{1\leq k\leq n}I(z_{k})$,
and $I(\mathcal{P}_{n})=\sup\limits_{p\in\mathcal{P}_{n}}I(p)$.\\
It was showed that there exists an extremal polynomial $p^{\ast}_{n}$,
 i.e., $I(\mathcal{P}_{n})=I(p^{\ast}_{n})$ and that $p^{\ast}_{n}$ has at least one zero on
each subarc of the unit circle of length $\pi$(see [1]).\\
It will suffice to prove the Sendov conjecture assuming $p$ is an extremal polynomial of the following form,
\bi
 p(z)=(z-a)\prod\limits_{k=1}^{n-1}(z-z_{k}),\ |z_{k}|\leq 1(k=1,\cdots,n-1)
 \ei
with\\
$p^{\prime}(z)=n\prod\limits_{j=1}^{n-1}(z-\zeta_{j}),\ \ |\zeta_{j}|\leq 1(j=1,\cdots,n-1),\ \ a\in [0,1]$.\\
Let $r_{k}=|a-z_{k}|,\ \rho_{j}=|a-\zeta_{j}|$ for
$k,j=1,2,\cdots,n-1$. By relabeling we suppose that \bi
\rho_{1}\leq\rho_{2}\leq \cdots\leq\rho_{n-1},\  r_{1}\leq r_{2}\leq\cdots\leq r_{n-1} . \ei
 We have(see [3])
 \bi
 2\rho_{1}\sin(\frac{\pi}{n})\leq r_{k}\leq 1+a,\ \
k=1,2,\cdots,n-1.
 \ei
 {\bf Sendov conjecture.} The disk $|z-a|\leq 1$ contains a zero of $p^{\prime}(z)$.\\
 If $a=0$  or $a=1$, Sendov conjecture is true(see [4]), we suppose $a\in (0,1)$.\\
 In this paper, we obtain the following theorem.\\
 {\bf Theorem.} For $n=9$, the disk $|z-a|\leq 1$ contains a zero of $p^{\prime}(z)$.\\
 Throughout the paper, we assume that $z_{k}\neq 0, \ \ z_{k}\neq a,\ k=1,\cdots,n-1$ and that $p(z)$ is extremal in the form (1.1): $I(\mathcal{P}_{n})=I(p)=I(a)=\rho_{1}$.\\
 In [5], we have proved the theorem for $ a<0.845$, it will suffice to prove the theorem in the case $ a\geq 0.845$.

\section{ Basic Lemmas }
\setcounter{equation}{0}

{\bf Lemma 2.1.} If
$\ \ 1-(1-|p(0)|)^{\frac{1}{n}}\leq \lambda\leq\sin(\frac{\pi}{n})$ and
$\lambda<a,\ \rho_{1}\geq 1$, then there exists a critical point
$\zeta_{0}=a+\rho_{0}e^{i\theta_{0}}$ such that $ Re \zeta_{0}\geq
\frac{1}{2}(a-\frac{\lambda(\lambda+2)}{a})$.\\
This is the Theorem 1 of [5].\\
{\bf Lemma 2.2.} For $\rho>0$, we have\\
$\prod\limits_{r_{k}\geq \rho}r_{k}\rho^{-1}\leq \prod\limits_{\rho_{j}\geq \rho}\rho_{j}\rho^{-1}
\prod\limits_{2\sin\frac{\pi k}{n}\geq 1}2\sin\frac{\pi k}{n}.$\\
This is the Theorem 3 of [5] by taking $m=1$.\\
{\bf Lemma 2.3.} If $c_{k}(k=1,\cdots,N),\ m,\ M, \ C$ are positive constants with\\
$m\leq c_{k}\leq M,\ \prod\limits_{k=1}^{N}c_{k}\geq C$ and $m^{N}\leq C\leq M^{N}$, then
$$\sum\limits_{k=1}^{N}\frac{1}{c_{k}^{2}}\leq\frac{N-v}{m^{2}}+\frac{v-1}{M^{2}}+\{\frac{m^{N-v}M^{v-1}}{C}\}^{2},$$
where $v=\min\{j\in \mathbb{Z}:M^{j}m^{N-j}\geq C\}.$\\
Proof. See Lemma 7.3.9 of [6] and Lemma B of [3].\\
{\bf Lemma 2.4.} We have
$$ \prod\limits_{k=1}^{n-1}r_{k}=n\prod\limits_{j=1}^{n-1}\rho_{j}.$$
Proof. See [3] .\\
Let $$\gamma_{j}=\frac{\zeta_{j}-a}{a\zeta_{j}-1}\  and \ w_{k}=\frac{z_{k}-a}{az_{k}-1}. $$
We take $n=9$, then
\bi
\prod\limits_{j=1}^{8}|\gamma_{j}|\leq \frac{\prod\limits_{k=1}^{8}|w_{k}|}{9-4a^{2}/(1+a^{2})-6a}.
\ei
{\bf Lemma 2.5.} If $\rho_{1}> 1,\ $ and $\ \zeta_{0}=a+\rho_{0}e^{i\theta_{0}}$  is the critical point in Lemma 2.1,  $\gamma_{0}=\frac{\zeta_{0}-a}{a\zeta_{0}-1}$, then
$$ |\gamma_{0}|>\frac{1}{\sqrt{1+\lambda(\lambda+2)-a^{2}\lambda(\lambda+2)}}.$$
Proof. See Lemma 3.4 of [5].\\
{\bf Lemma 2.6.} If $|\gamma_{j}|\leq\frac{1}{1+a-a^{2}}$, then $\rho_{j}\leq 1$.\\
Proof. See Lemma 1 of [2].\\
Write $B=B(R)=\frac{R+a}{aR+1}$.\\
{\bf Lemma 2.7.} If $ a\in [0.845,1]$, and there exists $|z_{k}|\leq 0.46$,
 then $I(a)=\rho_{1}\leq 1$.\\
Proof. If $I(a)>1$ and there exists $|z_{k}|\leq R$, then,
by (2.1) and Lemma 2.5, there exists some $\gamma_{j_{0}}$,\\
$ \frac{|\gamma_{j_{0}}|^{7}}{\sqrt{1+\lambda(\lambda+2)-a^{2}\lambda(\lambda+2)}}\leq\prod\limits_{j=1}^{8}|\gamma_{j}|
\leq \frac{B}{9-4a^{2}/(1+a^{2})-6a},$\\
hence\\
$|\gamma_{j_{0}}|^{7}\leq\frac{\sqrt{1+\lambda(\lambda+2)-a^{2}\lambda(\lambda+2)}}{9-4a^{2}/(1+a^{2})-6a}B$.\\
By Lemma 2.6, it suffices to show
\bi
 (9-4a^{2}/(1+a^{2})-6a)\frac{aR+1}{ R+a}- \sqrt{1+(1-a^{2})\lambda(\lambda+2)}(1+a-a^{2})^{7}\geq 0.
\ei
We consider the conditions for $\lambda$,\\
$|p(0)|=a\prod\limits_{k=1}^{8}|z_{k}|\leq aR,$
$\ \ 1-(1-|p(0)|)^{\frac{1}{9}}\leq \lambda\leq\sin(\frac{\pi}{9}),$\\
we choose
$$\lambda=1-(1-aR)^{\frac{1}{9}},$$
and $R$ satisfies $R\leq a^{-1}(1-(1-\sin(\frac{\pi}{9}))^{9}).$\\
Taking $R=0.46$, we obtain (2.2), the lemma follows.\\
We have\\
 $\frac{p^{\prime}(0)}{p(0)}=-(\frac{1}{a}+\sum\limits_{k=1}^{8}\frac{1}{z_{k}}),$  and
$9\prod\limits_{j=1}^{8}|\zeta_{j}|=(a\prod\limits_{k=1}^{8}|z_{k}|)|\frac{1}{a}+\sum\limits_{k=1}^{8}\frac{1}{z_{k}}|.$\\
Let $\Delta=Re(\frac{1}{a}+\sum\limits_{k=1}^{8}\frac{1}{z_{k}}),\ \sigma=\sum\limits_{k=1}^{8}\frac{1}{r_{k}^{2}}$, then\\
 $9\prod\limits_{j=1}^{8}|\zeta_{j}|\geq -\Delta a\prod\limits_{k=1}^{8}|z_{k}|$.\\
{\bf Lemma 2.8.} If $\rho_{1}>1,\ a\in [0.845,1]$, and $|z_{k}|\in [0.46, 1], \ k=1,\cdots,8$,  then
$$\Delta\leq -\frac{8}{a}+8a+\frac{9}{8a}(1-a^{2})\sigma.$$
Proof. See Lemma 3.9 of [5].\\
{\bf Lemma 2.9.} Let $m=\frac{1}{4},\ a\in [0.845,1], \ f(x)=\frac{x^{2}-1}{(1-x^{m})(a+x)^{2}}$, then
$f^{\prime}(x)>0$, for $x\in [0.46,1)$.\\
Proof. We have  $f^{\prime}(x)=\frac{Y}{(1-x^{m})^{2}(a+x)^{3}}$, where \\
$Y=((m-2)x^{m+1}-mx^{m-1}+2x)a+mx^{m+2}-(2+m)x^{m}+2$.\\
If $Y=0$, we will obtain $a>1$ for $x\in [0.46,1)$, hence $Y\neq 0$ for $a\in [0.845,1], x\in [0.46,1)$.
When $a=x=0.9$, $Y>0$, the lemma follows.\\
{\bf Lemma 2.10.} If $\rho_{1}\geq 1$, we have\\
$\prod\limits_{j=1}^{8}|\zeta_{j}|\leq (\prod\limits_{j=1}^{8}\rho_{j})(a^{2}-1+\frac{1}{4}\sum\limits_{k=1}^{8}\frac{|z_{k}|^{2}-a^{2}}{r_{k}^{2}})^{4}. $\\
Proof. See Lemma 3.11 of [5].\\
We will use the following conditions
\bi
\frac{x^{2}-1}{(1-x^{m})(a+x)^{2}}+(1-a^{2})(\sigma-4)\leq 0,\ \ x\in[0.46,1], a\in [0.845,1].
\ei
By Lemma 2.8, Lemma 2.9, and Lemma 2.10, we have the following lemma.\\
{\bf Lemma 2.11.} If $\sigma> 4,\ \rho_{1}>1$, condition (2.3) holds with $m=\frac{1}{4}$,\\
 and $|z_{k}|\in [0.46, 1], \ k=1,\cdots,8$, then
$$4^{4}(8-\frac{9}{8}\sigma)(1-a^{2})^{-3}\leq (\sigma-4)^{4}9\prod\limits_{j=1}^{8}\rho_{j}.$$
Proof. See Lemma 3.12 of [5], we use (2.3) to take the place of (3.6) there.\\
Write $R_{k}=r_{k}\prod\limits_{j=1}^{8}r_{j}^{-\frac{1}{8}},\ $
\bi
U^{\ast}(a)=(8-v^{\ast})(\frac{2\sin\frac{\pi}{9}}{1+a})^{-\frac{7}{4}}
+(v^{\ast}-1)(\frac{1+a}{9^{\frac{1}{8}}})^{-2}+\{(\frac{2\sin\frac{\pi}{9}}{1+a})^{\frac{7}{8}(8-v^{\ast})}
(\frac{1+a}{9^{\frac{1}{8}}})^{v^{\ast}-1} \}^{2},
\ei
where
\bi
v^{\ast}=\min\{j\in \mathbb{Z}: j\geq\{ 7\log(\frac{1+a}{2\sin\frac{\pi}{9}})\}
(\log\{\frac{(1+a)^{\frac{15}{8}}}{9^{\frac{1}{8}}(2\sin\frac{\pi}{9})^{\frac{7}{8}}}\})^{-1}\}.
\ei
{\bf Lemma 2.12.} If $\rho_{1}>1$,  then
$$ \sum\limits_{k=1}^{8}\frac{1}{R_{k}^{2}} \leq  U^{\ast}(a).$$
Proof. See Lemma 3.13 of [5].\\
Write
\bi
\sigma\leq U(a).
\ei
By Lemma 2.12, we obtain the following lemma.\\
{\bf Lemma 2.13.} If $\rho_{1}>1,\  4<\sigma\leq U(a)<\frac{64}{9},\ a\in [0.845,1]$, and $|z_{k}|\in [0.46, 1], \ k=1,\cdots,8$, then
$$\frac{4U(a)}{U(a)-4}(\frac{8-\frac{9}{8}U(a)}{(1-a^{2})^{3}})^{\frac{1}{4}}\leq U^{\ast}(a).$$
Proof. See Lemma 3.14 of [5], we use Lemma 2.11 to take the place of Lemma 3.12 there.
\section{ Proof of the Theorem }
\setcounter{equation}{0}
Suppose $n=9$, we want to give new upper bound to $U(a)$, this is the central part of the paper.\\
Write\\
$n_{1}=\#\{z_{k}: r_{k}<1\},\ n_{2}=\#\{z_{k}: r_{k}\geq 1\}, n_{1}+n_{2}=8$, and\\
$q=\prod\limits_{r_{k}<1}r_{k}$,
$\sigma=\sum\limits_{k=1}^{8}\frac{1}{r_{k}^{2}}=\sum\limits_{r_{k}<1}\frac{1}{r_{k}^{2}}+\sum\limits_{r_{k}\geq 1}\frac{1}{r_{k}^{2}}=\sum_{A}+\sum_{B}$, where\\
$\sum_{A}=\sum\limits_{r_{k}<1}\frac{1}{r_{k}^{2}}, \sum_{B}=\sum\limits_{r_{k}\geq 1}\frac{1}{r_{k}^{2}}$.\\
As we have pointed out in [5], by combining Lemma 2.2 and Lemma 2.3, new results may be obtained.\\
Suppose that $\rho_{1}>1$, by making use of Lemma 2.13 and the new upper bound of $U(a)$,
we will get a contradiction.\\
If $\rho_{1}>1$, taking $\rho=1$ in Lemma 2.2, we have\\
$
\prod\limits_{r_{k}\geq 1}r_{k}\leq \prod\limits_{j=1}^{8}\rho_{j}
\prod\limits_{2\sin\frac{\pi k}{9}\geq 1}2\sin\frac{\pi k}{9},
$\\
by Lemma 2.4,
$9\leq \prod\limits_{r_{k}<1}r_{k}\prod\limits_{2\sin\frac{\pi k}{9}\geq 1}2\sin\frac{\pi k}{9}$,
hence,
\bi
q\geq 9(\prod\limits_{2\sin\frac{\pi k}{9}\geq 1}2\sin\frac{\pi k}{9})^{-1}=(2\sin\frac{\pi}{9})^{2}.
\ei
By Lemma 2.4, in the sum $\sum_{B}$, we have
\bi
\prod\limits_{r_{k}\geq 1}r_{k}\geq 9q^{-1}.
\ei
In the sum $\sum_{A}$, we have
\bi
\prod\limits_{r_{k}<1}r_{k}=q<1.
\ei
By (1.2), (1.3), and Lemma 2.4, we obtain\\
 $\prod\limits_{k=1}^{8}r_{k}=9\prod\limits_{j=1}^{8}\rho_{j}>9,$\\
  then\\
 $2^{n_{2}}\geq (1+a)^{n_{2}}\geq\prod\limits_{r_{k}\geq 1}r_{k}\geq 9q^{-1}>9,$\\
 hence $n_{2}\geq 4$, and $n_{1}\leq 4$.\\
We will consider four subcases for $n_{1},\ a\in [0.845,1)$.\\
(i) $n_{1}=0, n_{2}=8$\\
$\sum_{A}=0=U_{A}$. say \\
In the sum $\sum_{B}$, we have $1\leq r_{k}\leq 1+a$, hence by Lemma 2.3\\
$\sum_{B}\leq n_{2}-v_{2}+\frac{v_{2}-1}{(1+a)^{2}}+\{\frac{(1+a)^{v_{2}-1}}{9}\}^{2}=U_{B},$ say\\
where $v_{2}=\min\{j\in\mathbb{Z}: (1+a)^{j}\geq 9\}$,\\
$U(a)=U_{A}+U_{B}$,\\
using the bound of $U(a)$, we get the desired contradiction to Lemma 2.13.\\
(ii) $n_{1}=1, n_{2}=7$\\
$\sum_{A}=q^{-2}=U_{A}$. say \\
In the sum $\sum_{B}$, we have $1\leq r_{k}\leq 1+a$, hence by Lemma 2.3\\
$\sum_{B}\leq n_{2}-v_{2}+\frac{v_{2}-1}{(1+a)^{2}}+\{\frac{(1+a)^{v_{2}-1}}{9q^{-1}}\}^{2}=U_{B},$ say\\
where $v_{2}=\min\{j\in\mathbb{Z}: (1+a)^{j}\geq 9q^{-1}\}$,\\
$U(a)=U_{A}+U_{B}$,\\
using the bound of $U(a), q\geq 2\sin\frac{\pi}{9}$, we get the desired contradiction to Lemma 2.13.\\
(iii) $n_{1}=2, n_{2}=6 $ or  $n_{1}=3, n_{2}=5 $\\
In the sum $\sum_{A}$, we have $2\sin\frac{\pi}{9}\leq r_{k}<1$, hence by Lemma 2.3\\
$\sum_{A}\leq
\frac{n_{1}-v_{1}}{(2\sin\frac{\pi}{9})^{2}}+v_{1}-1+\{\frac{(2\sin\frac{\pi}{9})^{n_{1}-v_{1}}}{q}\}^{2}=U_{A},$ say\\
where $v_{1}=\min\{j\in\mathbb{Z}: j\geq n_{1}-\frac{\log q}{\log(2\sin\frac{\pi}{9})}\}$.\\
In the sum $\sum_{B}$, we have $1\leq r_{k}\leq 1+a$, hence by Lemma 2.3\\
$\sum_{B}\leq n_{2}-v_{2}+\frac{v_{2}-1}{(1+a)^{2}}+\{\frac{(1+a)^{v_{2}-1}}{9q^{-1}}\}^{2}=U_{B},$ say\\
where $v_{2}=\min\{j\in\mathbb{Z}: (1+a)^{j}\geq 9q^{-1}\}$,\\
$U(a)=U_{A}+U_{B}$,\\
using the bound of $U(a), q\geq (2\sin\frac{\pi}{9})^{2}$, we get the desired contradiction to Lemma 2.13.\\
(iv) $n_{1}=4, n_{2}=4 $\\
In the sum $\sum_{A}$, we have $2\sin\frac{\pi}{9}\leq r_{k}<1$, hence by Lemma 2.3\\
$\sum_{A}\leq
\frac{n_{1}-v_{1}}{(2\sin\frac{\pi}{9})^{2}}+v_{1}-1+\{\frac{(2\sin\frac{\pi}{9})^{n_{1}-v_{1}}}{q}\}^{2}=U_{A},$ say\\
where $v_{1}=\min\{j\in\mathbb{Z}: j\geq n_{1}-\frac{\log q}{\log(2\sin\frac{\pi}{9})}\}$.\\
In the sum $\sum_{B},\ q(1+a)^{n_{2}}\geq 9$, we have $1\leq r_{k}\leq 1+a$, hence by Lemma 2.3\\
$\sum_{B}\leq n_{2}-v_{2}+\frac{v_{2}-1}{(1+a)^{2}}+\{\frac{(1+a)^{v_{2}-1}}{9q^{-1}}\}^{2}=U_{B},$ say\\
where $v_{2}=\min\{j\in\mathbb{Z}: (1+a)^{j}\geq 9q^{-1}\}$.\\
$U(a)=U_{A}+U_{B}$,\\
using the bound of $U(a), q\geq (2\sin\frac{\pi}{9})^{2}$, we get the desired contradiction to Lemma 2.13.\\
Finally, by Lemma 2.7, we obtain the theorem.\\
The method of this paper can be used to obtain new results for $n>9$, we leave this case to the reader.

{\small E-mail:\ mengzzh@126.com}
\end{document}